\documentclass[makeidx]{amsart}
\usepackage{amsmath}
\usepackage{amssymb}
\usepackage{epsfig}
\usepackage{hyperref}
\newtheorem{Proposition}{Proposition}
  \newtheorem{Remark}{Remark}
  \newtheorem{Corollary}[Proposition]{Corollary}

  \newtheorem{Theorem}{Theorem}

\newcommand {\z}{{\noindent}}
\def\blackslug{\hbox{\hskip 1pt \vrule width 4pt height 8pt depth 1.5pt
\hskip 1pt}}
\def\qed{\quad\blackslug\lower 8.5pt\null\par}

 \def\RR{\mathbb{R}}
 \def\NN{\mathbb{N}}
\def\ZZ{\mathbb{Z}}

\def\Re{\mathrm{Re}}

\def\bfu{\mathbf{u}}
\def\bfv{\mathbf{v}}
\def\bp{\boldsymbol{\phi}}

\def\pr#1#2{\frac{\partial{#1}}{\partial{#2}}}
\def\r{\rho}
\def\th{\theta}
\def\D1{\mathcal{D}}
\def\l{\lambda}
\def\Bi{\Big |}
\def\vfi{\varphi}

\def\g{\gamma}

\def\g{\gamma}
\def\r{\rho}
\def\dn{\l^2+(n/2-\g)^2}
\def\D{\mathcal{D}}
\def\t{\theta}
\def\bv{\mathbf{v}}
\def\bfu{\mathbf{u}}
\def\l{\lambda}
\def\an{\alpha_\nu}

\def\vf{{v_\vfi}}
\def\vr{{v_\r}}
\def\wf{{w_\vfi}}
\def\wr{{w_\r}}

\def\bfw{\mathbf{w}}

\makeindex
\begin{document}

\author{O. Costin$^1$ and V. Maz'ya$^2$ }\title{Sharp Hardy-Leray inequality
  for axisymmetric divergence-free fields} \gdef\shortauthors{O.  Costin \& V.
  Maz'ya} \gdef\shorttitle{Sharp Hardy-Leray inequality} \thanks{$1$.
  Department of
  Mathematics, Ohio State University,  231 West 18th Avenue, Columbus, OH 43210.\\
  \phantom{mm}$2.$ Department of Mathematics, Ohio State University, Columbus,
  OH 43210, USA;\ Department of Mathematical Sciences, University of
  Liverpool, Liverpool L69 3BX, UK;\ Department of Mathematics, Link\"oping
  University, Link\"oping, SE-581 83, Sweden}
 \maketitle
\date{}

\bigskip

\begin{abstract}
  We show that the sharp constant in the classical $n$-dimensional Hardy-Leray
  inequality can be improved for axisymmetric divergence-free fields, and find
  its optimal value. The same result is obtained for $n=2$ without the
  axisymmetry assumption.
\end{abstract}

\section{Introduction}

 \vskip 0.5cm \setcounter{section}{1} Let $\bf u$ denote a
 $C_0^{\infty}(\RR^n)$ vector field in $\RR^n$. The following $n$-dimensional
 generalization of the one-dimensional Hardy inequality \cite{Hardy}, $ $
 \vskip 0.05cm \begin{equation}
  \label{eq:eq1}
  \int_{\RR^n}\frac{|\bfu|^2}{|x|^{2}} dx\le \frac{4}{(n-2)^2} \int_{\RR^n}{|\nabla \bfu|^2}dx
\end{equation} $ $\vskip 0.05cm\z
appears for $n=3$ in the pioneering Leray's paper on the
Navier-Stokes equations \cite{Leray}. The constant factor on the right-hand side is sharp.
Since one frequently deals with divergence-free fields in hydrodynamics, it
is natural to ask whether this restriction can improve the constant in
(\ref{eq:eq1}).

We show in the present paper that this is the case indeed if  $n>2$ and the
vector field $\bf u$ is axisymmetric by proving that the aforementioned
constant can be replaced by the (smaller) optimal value $ $ \vskip 0.05cm
\begin{equation}
  \label{eq:con}
  \frac{4}{(n-2)^2}\left(1-\frac{8}{(n+2)^2}\right)
\end{equation} $ $\vskip 0.05cm\z
which, in particular, evaluates to $68/25$ in three dimensions. This result
is a special case of a more general one concerning a divergence-free
improvement  of the multi-dimensional
sharp Hardy inequality$ $ \vskip
  0.05cm
\begin{equation}
  \label{eq:g1}
  \int_{\RR^n}|x|^{2\g-2}|\bfu|^2 dx\le \frac{4}{(2\g+n-2)^2}  \int_{\RR^n}|x|^{2\g}|\nabla \bfu|^2 dx
\end{equation}$ $ \vskip
  0.05cm

  Let $\boldsymbol{\phi}$ be a point 
on the $(n-2)$-dimensional unit sphere $S^{n-2}$ with  spherical coordinates
$\{\theta_j\}_{j=1,...,n-3}$ and $\varphi$, where $\theta_j\in (0,\pi)$ and
$\vfi\in [0,2\pi)$. A point $x\in\RR^n$ is represented as a triple
$(\r,\th,\bp)$, where $\rho>0$ and $\theta\in [0,\pi]$.  Correspondingly, we
write $\bfu=(u_\rho,u_\theta,\bfu_{\bp})$ with
$\bfu_{\bp}=(u_{\theta_{n-3}},...,u_{\theta_1},u_\vfi)$. 

The {\em condition of axial
symmetry} means that $\bfu$ depends only on $\rho$ and $\theta$.

\bigskip

\z For higher dimensions, our result is as follows. 
\vskip 0.3cm 

\begin{Theorem}\label{T1}
  Let  $\g\ne 1-n/2$, $n>2$, and
  let $\bfu$ be an axisymmetric divergence-free vector field in $C_0^{\infty}
  (\RR^n)$. We assume that $\bfu(\mathbf{0})=\bf 0$ for $\g<1-n/2$.  Then 
 \begin{equation}
  \label{eq:eq2}
  \int_{\RR^n}|x|^{2\g-2}|\bfu|^2dx\le C_{n,\g}\int_{\RR^n}|x|^{2\g}|{\nabla \bfu|^2}dx
\end{equation} $ $\vskip 0.05cm\z with
the best value of
  $C_{n,\g}$  given by
\begin{equation}
  \label{eq:C}
  C_{n,\g}=\frac{4}{(2\g+n-2)^2}\left(1-\frac{2}{n+1+(\g-n/2)^2}\right)
\end{equation}
for $\gamma\le 1$, and by
\begin{multline}
  \label{eq:ref2}
 {C^{-1}_{n,\g}}=\left(\frac{n}{2}+\g-1\right)^2\\+\min\Big\{ n-1,\ 2+\min_{x\ge 0}\Big(x+\frac{4(n-1)(\g-1)}{x+n-1+(\g-n/2)^2}\Big)\Big\}
\end{multline}
\bigskip
for $\gamma>1$.
\end{Theorem}

\vskip -0.6cm

The two minima in  (\ref{eq:ref2}) can be calculated in closed form, but their
expressions for arbitrary dimensions turn out to be  unwieldy, and we
omit them. 

However, the formula for $C_{3,\g}$ is  simple.
\begin{Corollary}
For $n=3$ inequality {\rm (\ref{eq:eq2})} holds with the best constant
\begin{equation}
  \label{newform}
  C_{3,\g}=\left\{\begin{array}{lll}
\displaystyle \frac{4}{(2\g+1)^2}\cdot \frac{2+\left(\g-3/2\right)^2}{4+(\g-3/2)^2},
 \ \ \text{for $\gamma\le 1$}\\ \\   \displaystyle \frac{4}{8+(1+2\g)^2},\ \ \ \ \ \ \ \ \text{for $\gamma>1$.}\end{array}\right.
\end{equation}
\vskip 0.05cm 
\end{Corollary}
For $n=2$, we obtain the sharp constant in (\ref{eq:eq2}) without
axial symmetry  of the vector field.

\begin{Theorem}\label{T2}
  Let $\gamma\ne 0$, $n=2$, and let $\bf u$ be a divergence-free vector field
  in $C_0^{\infty}(\RR^2)$.  We assume that $\bfu(\mathbf{0})=\bf 0$ for
  $\g<0$.  Then inequality {\rm (\ref{eq:eq2})} holds with the best constant
\begin{equation}
  \label{eq:C1}
  C_{2,\g}=\left\{\begin{array}{lll}\displaystyle \gamma^{-2}
      \frac{1+(1-\g)^2}{3+(1-\g)^2}& \text{for } \gamma\in [-\sqrt{3}-1,\sqrt{3}-1]
      \\ \\ (\gamma^2+1)^{-1}& \text{otherwise. } \end{array}\right.
\end{equation}
\end{Theorem}

\section{Proof of Theorem 1}
In the spherical coordinates introduced above, we have
$ $ \vskip 0.05cm \begin{multline}
  \label{eq:140}
  {\rm div}\,\, \mathbf{u}=\r^{1-n}\pr{}{\r}\left(\r^{n-1}u_\r\right)+\r^{-1}(\sin\th)^{2-n}\pr{}{\th}\left((\sin\th)^{n-2}u_\th\right)\\+\sum_{k=1}^{n-3}(\r\sin\th\sin\th_{n-3}\cdots\sin\th_{k+1})^{-1}(\sin\th_k)^{-k}\pr{}{\th_k}\left((\sin\th_k)^ku_{\th_k}\right)\\+(\r\sin\th\sin\th_{n-3}\cdots\sin\th_1)^{-1}\pr{u_\vfi}{\vfi}
\end{multline}$ $\vskip 0.05cm\z
Since the components $u_\vfi$ and $u_{\th_k}$, $k=1,...,n-3$, depend only on
$\r$ and $\th$, (\ref{eq:140}) becomes
\begin{multline}
  \label{eq:141}
  {\rm div}\,\, \mathbf{u}=\r^{1-n}\pr{}{\r}\left(\r^{n-1}u_\r(\r,\th)\right)+\r^{-1}(\sin\th)^{2-n}\pr{}{\th}\left((\sin\th)^{n-2}u_\th(\r,\th)\right)\\+\sum_{k=1}^{n-3}k(\sin\th_{n-3}\cdots\sin\th_{k+1})^{-1}\cot\th_k\frac{u_{\th_k}(\r,\th)}{\r\sin\th}
\end{multline}$ $\vskip 0.05cm\z By the linear independence of the functions

$$1,\,\,\,(\sin\th_{n-3}\cdots\sin\th_{k+1})^{-1}\cot\th_k,\ \ k=1,...,n-3$$ 
\vskip 0.2cm

\z the
divergence-free condition is equivalent to the collection of $n-2$ identities
$ $ \begin{equation}
  \label{eq:142}
  \r\pr{u_\r}{\r}+(n-1)u_\r+\left(\pr{}{\th}+(n-2) \cot\th\right)u_\th=0
\end{equation} $ $\vskip 0.05cm\z
\begin{equation}
  \label{eq:143}
  u_{\th_k}=0,\ \ \ k=1,...,n-3
\end{equation} $ $ \vskip 0.05cm \z

If the right-hand side of (\ref{eq:eq2}) diverges, there is nothing to prove.
Otherwise, the matrix $\nabla\mathbf{u}$ is $O(|x|^{m})$, with $m>-\g-n/2$, as
$x\to 0$. Since $\mathbf{u}(\mathbf{0})=\mathbf{0}$, we have
$\mathbf{u}(x)=O(|x|^{m+1})$ ensuring the convergence of the integral on the
left-hand side of (\ref{eq:eq2}). We introduce the vector field $ $ \vskip
0.05cm \begin{equation}
  \label{eq:eqv}
  \bfv (x)={\bfu (x)}|x|^{\g-1+n/2}
\end{equation} $ $\vskip 0.05cm\z
The inequality (\ref{eq:eq2}) becomes 
$ $ \vskip 0.05cm \begin{equation}
  \label{eq:eq4}
\left(\frac{1} {C_{n,\g}}-\left(\frac{n}{2}+\g-1\right)^2\right)\int_{\RR^n}\frac{|\bfv|^2}{|x|^n} dx\le\int_{\RR^n}\frac{|\nabla \bfv|^2}{|x|^{n-2}}dx
\end{equation} $ $\vskip 0.05cm\z
The condition ${\rm div\,\,} \bf u\rm=0$ is equivalent to
 \begin{equation}
  \label{eq:eqdiv}
  \rho\,\,\text{div}\,\,\bfv =\left(\frac{n-2}{2}+\g\right){v_\rho}
\end{equation}

To simplify the exposition, we assume first that $\bf v_{\bp}=0$. Now,
(\ref{eq:eqdiv}) can be written as
$ $ \vskip 0.05cm \begin{equation}
  \label{eq:eqdiv1}
  \rho\frac{\partial v_\rho}{\partial{\rho}}+\left(\frac{n}{2}-\g\right)v_\rho
  +\mathcal{D}v_\theta=0
\end{equation} 
where
\begin{equation}
  \label{eq:eqD}
  \mathcal{D}:=\frac{\partial }{\partial \theta}+ (n-2)\cot \theta\,\,
\end{equation} $ $\vskip 0.05cm\z
Note that $\mathcal{D}$ is the adjoint of $-\partial/\partial\theta$ with
respect
to the scalar product
$$\int_0^{\pi}f(\theta)\overline{g(\theta)}(\sin\theta)^{n-2}d\theta$$
A straightforward though lengthy calculation yields
$ $\vskip 0.02cm \begin{multline}
  \label{eq:nab}
  \rho^2 |\nabla \mathbf{v}|^2= \rho^2 \Big(\pr{v_\r}{\r}\Big)^2 + \rho^2 \Big(\pr{v_\th}{\r}\Big)^2
+\Big(\pr{v_\r}{\th}\Big)^2+\Big(\pr{v_\th}{\th}\Big)^2
\\+v_\th^2+(n-1)v_\r^2 +(n-2)(\cot\th)^2 v_\th^2+2\left(v_\r\D  v_\th-v_\th\pr{v_\r}{\th}\right)
\end{multline}$ $\vskip 0.05cm\z
Hence
\begin{multline}\label{210}
   \rho^2 \int_{S^{n-1}} |\nabla \mathbf{v}|^2ds=
 \int_{S^{n-1}}\Big\{ \rho^2 \Big( \pr{v_\r}{\r}\Big)^2 + 
\Big(\pr{v_\th}{\th}\Big)^2+\rho^2 \Big(\pr{v_\th}{\r}\Big)^2+
\Big(\pr{v_\r}{\th}\Big)^2\\ \\ +v_\th^2+(n-1)v_\r^2 +(n-2)(\cot\th)^2 v_\th^2+4 v_\r\D  v_\th\Big\}ds
\end{multline}$ $\vskip 0.02cm\z
Changing the  variable $\rho$  to $t=\log\rho$, and applying the  Fourier
transform with respect to
$t$,
 $$\mathbf v(t,\theta)\mapsto{\mathbf w}(\lambda,\theta)$$ we derive
\begin{multline}\label{111}
  \int_{\RR^n}
  \frac{|\nabla\mathbf{v}|^2}{|x|^{n-2}}dx\\=\int_\RR\int_{S^{n-1}}\Big\{(\l^2+n-1)|w_\r|^2+(\l^2-n+3)|w_\th|^2\\+\Bi\pr{w_\r}{\th}\Bi^2+\Bi\pr{w_\th}{\th}\Bi^2+(n-2)(\sin\th)^{-2}|w_\th|^2+4\Re( \overline{w}_\r\D  {w_\th})\Big\}dsd\l
\end{multline}$ $\vskip 0.05cm\z
and
\begin{equation}
  \label{eq:Id2}
  \int_{\RR^n}\frac{|\mathbf{v}|^2}{|x|^n}dx=\int_{\RR}\int_{S^{n-1}}|\mathbf{w}|^2 ds d\lambda
\end{equation}
From (\ref{eq:eqdiv}), we obtain
$ $ \vskip 0.05cm \begin{equation}
  \label{eq:il}
 w_\rho= -\frac{\mathcal{D}w_\theta}{i\lambda+n/2-\g}
\end{equation}
which implies
 \begin{equation}
  \label{eq:124}
  |w_\r|^2=\frac{|\D  w_\th|^2}{\dn}
\end{equation}
and 
 \begin{equation}
  \label{eq:125}
  \Re( \overline{w}_\r\D  {w_\th})=-\frac{(n/2-\g)|\D  w_\th|^2}{\dn}
\end{equation} $ $\vskip 0.05cm\z
Introducing this into (\ref{111}), we arrive at the identity
\begin{multline}
  \label{215}
  \int_{\RR^n}\frac{|\nabla\bv |^2}{|x|^{n-2}}dx=\int_0^{\infty}\int_0^{\pi}
  \Bigg\{(\l^2+n-1)\frac{|\D
    w_\t|^2}{\l^2+(n/2-\g)^2}+(\l^2-n+3)|w_\t|^2\\+\left|\pr{w_\t}{\t}\right|^2+(n-2)(\sin\t)^{-2}|w_\t|^2+\frac{1}{\l^2+(n/2-\g)^2}\left|\pr{}{\t}\D
  w_\t\right|^2
  \\-4\left(\frac{n}{2}-\g\right)\frac{|\D w_\t|^2}{\l^2+(n/2-\g)^2}\Bigg\}d\t d\l
\end{multline}
We simplify the right-hand side of (\ref{215}) to obtain
\begin{multline}
  \label{mul1}
\int_{\RR^n}\frac{|\nabla\bv |^2}{|x|^{n-2}}dx=
\int_0^{\infty}\int_0^{\pi}\Bigg\{\left(\frac{-n-1+\l^2+4\g}{\l^2+(n/2-\g)^2}+1\right)|\D
w_\t|^2\\+(\l^2-n+3)|w_\t|^2+\frac{1}{\l^2+(n/2-\g)^2}\left|\pr{}{\t}\D
  w_\t\right|^2\Bigg\} d\t d\l
\end{multline}
$ $ \vskip 0.05cm 
On the other hand, by (\ref{eq:Id2}) and (\ref{eq:il}) 
$ $ \vskip 0.05cm \begin{equation}
  \label{eq:g14}
  \int_{\RR^n}\frac{|\bv
    |^2}{|x|^{n-2}}dx=\int_0^{\infty}\int_0^{\pi}\left(\frac{|\D
    w_\t|^2}{\l^2+(n/2-\g)^2}+|w_\t|^2\right)d\t d\l
\end{equation}$ $ \vskip 0.05cm 
Defining the self-adjoint operator 
$ $ \vskip 0.05cm \begin{equation}
  \label{eq:128}
  T:=-\pr{}{\th}\D
\end{equation} $ $\vskip 0.05cm\z
or, equivalently,
\begin{equation}
  \label{eq:214}
  T=-\delta_\th+\frac{n-2}{(\sin\th)^2}
\end{equation}$ $ \vskip 0.05cm 
\z where $\delta_\th$ is the $\th$-part of the Laplace-Beltrami operator on
$S^{n-1}$, we write (\ref{mul1}) and (\ref{eq:g14}) as $ $ \vskip 0.05cm
\begin{equation}
  \label{1263}
\int_{\RR^n}
  \frac{|\nabla\mathbf{v}|^2}{|x|^{n-2}}dx=\int_\RR\int_{S^{n-1}}Q(\l,w_\th)dsd\l
\end{equation} $ $\vskip 0.05cm\z
and 
$ $ \vskip 0.05cm \begin{equation}
  \label{eq:146}
\int_{\RR^n}
  \frac{|\mathbf{v}|^2}{|x|^{n}}dx=\int_\RR\int_{S^{n-1}}q(\l,w_\th)dsd\lambda
\end{equation}$ $ \vskip 0.05cm 
\z respectively,  where $Q$ and $q$ are  sesquilinear  forms in $w_\t$, defined by
$ $ \vskip 0.05cm \begin{multline}
  \label{eq:145}
  Q(\l, w_\th)=
\left(\frac{-n-1+\l^2+4\g}{\dn}+1\right)T w_\th\cdot\overline{w_\th}\\+(\l^2-n+3)|w_\th|^2
+\frac{1}{\dn}|Tw_\th|^2
\end{multline} $ $\vskip 0.05cm\z
and
$ $ \vskip 0.05cm \begin{equation}
  \label{eq:147}
 q(\l,w_\th)= \frac{T w_\th\cdot\overline{w_\th}}{\dn}+|w_\th|^2 
\end{equation} $ $\vskip 0.05cm\z
The
eigenvalues of $T$ are $\alpha_\nu=\nu(\nu+n-2)$, $\nu\in\ZZ^+$. Representing $w_{\theta}$ 
as an expansion in eigenfunctions of $T$, we find
$ $ \vskip 0.05cm \begin{multline}
  \label{eq:1328}
  \inf_{w_\th}\,\,\,\frac{\displaystyle \int_\RR\int_{S^{n-1}}Q(\l,w_\th)dsd\l}{\displaystyle \int_\RR\int_{S^{n-1}}q(\l,w_\th)dsd\lambda}\\=\inf_{\lambda\in\RR}\inf_{\nu\in\NN}\frac{\displaystyle\left(\frac{-n-1+\l^2+4\g}{\displaystyle\dn}+1\right)\alpha_\nu+\l^2-n+3
+\frac{\displaystyle \alpha_\nu^2}{\dn}}{_{\frac{\displaystyle \alpha_\nu}{_{\displaystyle \dn+1}} }}
\end{multline} $ $\vskip 0.05cm\z
Thus our minimization problem reduces to finding
\begin{equation}
  \label{eq:e43}
  \inf_{x\ge 0} \inf_{\nu\in\NN} f(x,\alpha_\nu,\g)
\end{equation}
where 
\begin{equation}
  \label{eq:eq46}
  f(x,\alpha_\nu,\gamma)= x-n+3+\alpha_\nu\left(1-\frac{16(1-\gamma)}{4x+4\alpha_\nu+(n-2\g)^2}\right)
\end{equation}$ $ \vskip 0.05cm 
Since $\gamma\le 1$, it is clear that $f$ is increasing in
$x$, so the value (\ref{eq:e43}) is equal to 
$ $ \vskip 0.05cm \begin{equation}
  \label{eq:1234}
  \inf_{\nu\in\NN} f(0,\alpha_\nu,\g)= \inf_{\nu\in\NN} \left(3-n+\alpha_\nu\left(1-\frac{16(1-\g)}{4\alpha_\nu+(n-2\g)^2}\right)\right)
\end{equation}$ $ \vskip 0.05cm 
We have$ $ \vskip 0.05cm 
\begin{equation}
  \label{eq:eqder}
  \pr{}{\alpha_\nu}f(0,\an,\g)=1-\frac{16 (1-\g)(n-2\g)}{(4\an+(n-2\g)^2)^2}
\end{equation}$ $ \vskip 0.05cm \z 
Noting that $ $ \vskip 0.05cm 
\begin{equation}
  \label{eq:csch}
  4\an+(n-2\g)^2\ge 4(n-1)+(n-2\g)^2\ge 4\sqrt{n-1}(n-2\g)
\end{equation}$ $ \vskip 0.05cm 
\z we see that 
\begin{equation}
  \label{eq:ineq}
   \pr{}{\alpha_\nu}f(0,\an,\g)\ge 1-\frac{1-\g}{(n-1)(n-2\g)}>0
\end{equation} Thus the minimum of $f(0,\an,\g)$ is attained at $\alpha_1=n-1$ and equals
$ $ \vskip 0.05cm \begin{equation}
  \label{eq:frmf}
  3-n+(n-1)\left(1-\frac{16(1-\g)}{4(n-1)+(n-2\g)^2}\right)=\frac{2(\g-1+n/2)^2}{n-1+(\g-n/2)^2}
\end{equation}$ $ \vskip 0.05cm 
\z This completes the proof for the case $\bf v_{\bp}=0$. 

\vskip 0.3cm 

If we drop the assumption $\bf v_{\bp}=0$, then, to the integrand on the right-hand side of
(\ref{210}), we should add   the terms
$ $ \vskip 0.05cm \begin{equation}
  \label{eq:1500}
  \r^2\left(\pr{v_\vfi}{\r}\right)^2+\left(\pr{v_\vfi}{\th}\right)^2+
\left(\sin\th\sin\th_{n-3}\cdots\sin\th_1\right)^{-2}
v_\vfi^2
\end{equation} $ $\vskip 0.05cm
\z The expression in (\ref{eq:1500}) equals
$ $ \vskip 0.05cm 
\begin{equation}
  \label{eq:150}
 \r^2\left|\nabla(v_\vfi e^{i\vfi})\right|^2
\end{equation} $ $\vskip 0.05cm\z 
As a result, the right-hand side of (\ref{1263}) is augmented by 
$ $ \vskip 0.05cm \begin{equation}
  \label{eq:160}
  \int_\RR\int_{S^{n-1}}R(\lambda,w_\vfi)dsd\l
\end{equation} $ $\vskip 0.05cm\z
where
 \begin{equation}
  \label{eq:170}
  R(\l,w_\vfi)=\l^2|w_\vfi|^2+|\nabla_\omega(w_\vfi e^{i\vfi})|^2
\end{equation} $ $\vskip 0.05cm\z
with $\omega=(\th,\th_{n-3},...,\vfi)$. Hence,
$ $ \vskip 0.05cm \begin{equation}
  \label{eq:180}
 \inf_{\bf v}\,\,\, \frac{\displaystyle \int_{\RR^n}\frac{|\nabla\mathbf{v}|^2}{|x|^{n-2}}dx}{\displaystyle\int_{\RR^n}
  \frac{ |\mathbf{v}|^2}{|x|^{n}}dx}\ =\ \inf_{w_{\th},w_\vfi}\frac{\displaystyle \int_\RR\int_{S^{n-1}}\Big(Q(\l,w_\th)+R(\lambda,w_{\vfi})\Big)dsd\l}{\displaystyle \int_\RR\int_{S^{n-1}}\Big(q(\lambda,w_\th)+|w_\vfi|^2\Big)dsd\l}
\end{equation} $ $\vskip 0.05cm\z
Using the fact that  $w_\th$ and $w_\vfi$ are independent, the right-hand side
is the lesser of  (\ref{eq:1328})
and
$ $ \vskip 0.05cm \begin{equation}
  \label{eq:2100}
 \inf_{w_{\vfi}}\,\,\,\frac{\displaystyle \int_\RR\int_{S^{n-1}}R(\lambda,w_\vfi)dsd\l}{\displaystyle \int_\RR\int_{S^{n-1}}|w_\vfi|^2dsd\l}
\end{equation}$ $ \vskip 0.05cm 
\z  Since $w_\vfi e^{i\vfi}$ is orthogonal to one on
$S^{n-1}$, we have
$ $ \vskip 0.05cm \begin{equation}
  \label{eq:191}
  \int_{S^{n-1}}\left|\nabla_\omega\left(w_\vfi e^{i\vfi}\right)\right|^2ds
\ge (n-1)\int_{S^{n-1}}|w_\vfi|^2ds
\end{equation} $ $\vskip 0.05cm\z
Hence the infimum in (\ref{eq:2100})  is at most $n-1$, which exceeds the value in
(\ref{eq:frmf}). The result follows for $\gamma \le 1$. 

For $\gamma>1$ the proof is similar. 
Differentiation of $f$ in $\alpha_\nu$ gives $ $ \vskip 0.05cm
\begin{equation}
  \label{eq:g34}
 1+ \frac{16(\g-1)((n-2\g)^2+4x)}{(4x+4\alpha_\nu+(n-2\g)^2)^2}
\end{equation}$ $ \vskip 0.05cm 
\z which is positive. Hence the role of the value (\ref{eq:frmf}) is played by
the smallest value of $f(\cdot,n-1,\gamma)$ on $\RR^+$. Therefore, 
$ $ \vskip 0.05cm\begin{equation}
  \label{eq:eqnew}
   \inf_{\bf v}\,\,\, \frac{\displaystyle \int_{\RR^n}\frac{|\nabla\mathbf{v}|^2}{|x|^{n-2}}dx}{\displaystyle\int_{\RR^n}
  \frac{ |\mathbf{v}|^2}{|x|^{n}}dx}\ =2+\min_{x\ge 0}\Big(x+\frac{4(n-1)(\g-1)}{x+n-1+(\g-n/2)^2}\Big)
\end{equation}
$ $ \vskip 0.1cm \z  The proof is complete.

\bigskip

{\em Proof of  Corollary 1}. We need to consider only  $\g>1$. It follows
directly from (\ref{eq:ref2})  that
$${C^{-1}_{3,\g}}=\left(\frac{3}{2}+\g-1\right)^2+2$$
which gives the result.
\z \begin{Remark}\label{R1} {\rm Using (\ref{eq:il}), we see that a minimizing
    sequence $\{\mathbf{v}_k\}_{k\ge 1}$ which shows the sharpness of
    inequality (\ref{eq:eq2}) with the constant
    (\ref{eq:C}) can be obtained by taking
    $\mathbf{v}_k=(v_{\r,k},v_{\th,k},\mathbf{0})$ with the Fourier transform
    ${\bf w}_k=(w_{\r,k},w_{\th,k},\mathbf{0})$ chosen as follows: $ $ \vskip
    0.05cm \begin{equation}
  \label{eq:1355}
  w_{\th,k}(\lambda,\th)=h_k(\lambda)\sin\theta,\ \
  w_{\r,k}(\lambda,\th)=\frac{1-n}{i\l+n/2-\g}h_k(\l)\cos\th
\end{equation} $ $\vskip 0.05cm\z
The sequence $\{|h_k|^2\}_{k\ge 1}$ converges in distributions to the
delta function at $\lambda=0$. The minimizing sequence that gives the value
(\ref{newform}) of $C_{3,\g}$ is

$$w_{\t,k}(\l,0)=0,\ \ w_{\rho,k}(\l,\theta)=0,\  \text{and}\ \ w_{\phi,k}(\l,\theta)=h_k(\l)
\sin\theta$$
$ $ \vskip 0.05cm 
\z where  $\{|h_k|^2\}_{k\ge 1}$ is as above.

\bigskip }
\end{Remark}

\section{Proof of Theorem 2.} 

The calculations are similar but simpler than those in the previous section. 
We start with the substitution $\bfv(x)=\bfu(x)|x|^{2\g}$ and write
(\ref{eq:eq4}) in the form
$ $ \vskip 0.05cm 
\begin{equation}
  \label{eq:c2}
  \frac{1}{C_{2,\g}}=\g^2+\inf_{\bfv}\frac{^{\displaystyle\int_{\RR^2} |\nabla
      \bfv|^2 dx}}{_{\displaystyle \int_{\RR^2}|\bfv |^2|x|^{-2}dx}}
\end{equation}
 $ $ \vskip 0.05cm \z 
In  polar coordinates $\rho$ and $\vfi$, with $\vfi\in [0,2\pi)$, we have$ $ \vskip 0.05cm 
\begin{multline}
  \label{eq:eqv2}
  \int_{\RR^2} |\nabla \bfv|^2 dx=\int_{\RR^2}\Big\{ |\nabla v_\r|^2+ |\nabla
  v_\vfi|^2+\rho^{-2}\left({\vr^2+\vf^2}-4\vr(\partial_\vfi v_\vfi)\right)
  \Big\}dx
\end{multline}$ $ \vskip 0.05cm \z 
Changing the  variable $\rho$  to $t=\log\rho$, and applying the  Fourier
transform $\bfv(\r,\vfi)\to\bfw(\l,\vfi)$, we obtain
$ $ \vskip 0.05cm \begin{equation}
  \int_{\RR}\int_0^{2\pi}\Big\{(\l^2+1)(|w_\r|^2+|w_\vfi|^2)+|\partial_\vfi
  w_\vfi|^2+|\partial_\vfi w_\r|^2-4\,(\partial_\vfi w_\vfi)\overline{\wr}\Big \}d\vfi d\l
\end{equation}$ $ \vskip 0.05cm \z
The divergence-free condition for $u$ becomes
$ $ \vskip 0.05cm \begin{equation}
  \label{eq:divf}
  w_\rho=-\frac{\partial_\vfi\wf}{i\l+1-\g}
\end{equation}$ $ \vskip 0.05cm 
\z
which yields
$ $ \vskip 0.05cm \begin{multline}
  \label{eq:eqvf2}
  \int_{\RR^2} |\nabla
      \bfv|^2 dx=\int_{\RR}\int_0^{2\pi}\Big\{ \Big(\frac{\l^2+4\g
          -3}{\l^2+(1-\g)^2}+1\Big)\left|\partial_\vfi
          \wf\right|^2\\+\frac{|\partial^2_\vfi w_\vfi|^2}{\l^2+(1-\g)^2}
+(\l^2+1)| w_\vfi|^2 \Big\}d\vfi d\l
\end{multline}
$ $ \vskip 0.05cm \z 
Analogously,
$ $ \vskip 0.05cm\begin{multline}
 \int_{\RR^2}|\bfv|^2
 |x|^{-2}dx=\int_{\RR}\int_0^{2\pi}(|w_\r|^2+|w_\vfi|^2)d\vfi d\l\\=
 \int_{\RR}\int_0^{2\pi}\Big(\frac{|\partial_\vfi
   w_\vfi|^2}{\l^2+(1-\g)^2}+|w_\vfi|^2\Big)d\vfi d\l
\end{multline}
$ $ \vskip 0.05cm \z
 Therefore, by (\ref{eq:c2}) 
 \begin{equation}
   \label{eq:newC}
   \frac{1}{C_{2,\g}}=\g^2+\inf_{x\ge 0}\,\,\inf_{\nu\in \NN \cup 0}f(x,\nu,\g)
 \end{equation}
where $ $ \vskip 0.05cm \z
\begin{equation}
  \label{eq:fff}
  f(x,\nu,\g)=x+1+\nu\left(1-\frac{4(1-\g)}{x+\nu+(1-\g)^2}\right)
\end{equation}$ $ \vskip 0.05cm \z
Let first $\g\le 1$. Then $f$ is  increasing in $x$, which implies
$f(x,\nu,\g)\ge f(0,\nu,\g)$. Since the derivative $ $ \vskip 0.05cm \z
\begin{equation}
  \label{eq:diff}
  \frac{\partial}{\partial\nu}f(0,\nu,\g)=1-\frac{4(1-\g)^3}{(\nu+(1-\g)^2)^2}
\end{equation}$ $ \vskip 0.05cm \z
is positive for $\nu\ge 2$, we need to compare only the values $f(0,0,\g)$,
$f(0,1,\g)$ and $f(0,2,\g)$. An elementary calculation shows that  both
$f(0,0,\g)$ and $f(0,2,\g)$ exceed 
$f(0,1,\g)$  for $\g\not\in
(-1-\sqrt{3},-1+\sqrt{3})$.

Let now $\g>1$. We have $ $ \vskip 0.05cm \z
\begin{equation}
  \label{eq:dfv}
  \frac{\partial}{\partial \nu} f(x,\nu,\g)=1+\frac{4(\g-1)(x+(1-\g)^2)}{(x+\nu+(1-\g^2))^2}>0
\end{equation} $ $ \vskip 0.05cm \z
 and therefore $f(x,\nu,\g)\ge f(x,0,\g)=x+1\ge 1$. The proof of Theorem 2 is
 complete. 
 \begin{Remark}{\rm 
   Minimizing sequences which give  $C_{2,\g}$ in (\ref{eq:C1}) can be
   chosen as follows:

$$ w_{\rho,k}(\l,\vfi)=0,\ \ \   \ \ \ w_{\vfi,k}(\l,\vfi)=h_k(\l)
$$
$ $ \vskip 0.05cm 
\z  for $\g\not\in
(-1-\sqrt{3},-1+\sqrt{3})$, and 
$ $ \vskip 0.05cm \begin{equation*}
  \label{eq:divf1}
  w_{\rho,k}=\frac{h_k(\l)\sin (\vfi-\vfi_0)}{i\l+1-\g},\ \ \  \ \ \  w_{\vfi,k}=
h_k(\l)\cos (\vfi-\vfi_0)
\end{equation*}$ $ \vskip 0.05cm

\z when $\g\in (-1-\sqrt{3},-1+\sqrt{3})$, for any constant $\vfi_0$. Here
$\{|h_k|^2\}_{k\ge 1}$ converges in distributions to the delta function at
$0$.}
 \end{Remark}
 \begin{Corollary}
    Let $\gamma\ne 0$. Denote by $\psi$ a real-valued scalar function in
     $C_0^{\infty}(\RR^2)$ and assume in addition that $\nabla
     \psi(\mathbf{0})=\mathbf{0}$ if $\gamma<0$. Then the sharp
     inequality
$ $ \vskip 0.05cm \begin{equation}
  \label{eq:divf123}
  \int_{\RR^2}|\nabla \psi|^2|x|^{2(\g-1)}dx\le C_{2,\g}\int_{\RR^2}\left(\psi_{x_1x_1}^2+2\psi_{x_1x_2}^2+\psi_{x_2x_2}^2\right)|x|^{2\g}dx
\end{equation}$ $ \vskip 0.05cm 
\z holds with $C_{2,\g}$ given in {\rm (\ref{eq:C1})}. 
 \end{Corollary}

 Indeed, for $n=2$, inequality (\ref{eq:eq2}) becomes (\ref{eq:divf123}) if
 $\psi$ is interpreted as a stream function of the vector field $\bf u$, {\em
   i.e.} $\mathbf{u}=\nabla\times \psi$.

\bigskip

\bigskip

\z {\bf Acknowledgments.} The authors are grateful to S. Tanveer and G. Luo for
providing some useful formulas. The work was partially supported
by the NSF grants DMS 0406193 and DMS 0600369 (O.C.), and DMS 0500029 (V.M.).

\end{document}